\documentclass[11pt,a4paper]{amsart}
\usepackage{amssymb}
\usepackage{mathrsfs} 
 \usepackage{amsmath}
\usepackage{amsfonts}
\usepackage[english]{babel}
\usepackage[T1]{fontenc}
\usepackage[latin1]{inputenc}
\usepackage{amsthm}
\makeatletter
\@namedef{subjclassname@2020}{%
  \textup{2020} Mathematics Subject Classification}
\makeatother
\newtheorem{theorem}{Theorem}[section]
\newtheorem{lemma}[theorem]{Lemma}

\theoremstyle{definition}

\newtheorem{prop}[theorem]{Proposition}
\theoremstyle{remark}
\newtheorem{remark}[theorem]{Remark}
\newtheorem{remarks}[theorem]{Remarks}
\numberwithin{equation}{section}
\usepackage{geometry}
\begin{document}

\title{Exotic closed subideals of algebras of bounded  operators}
 
\author[Hans-Olav Tylli]{Hans-Olav Tylli}
\address{Tylli: Department of Mathematics and Statistics, Box 68,
FI-00014 University of Helsinki, Finland}
\curraddr{}
\email{hans-olav.tylli@helsinki.fi}
\thanks{}

\author[Henrik Wirzenius]{Henrik Wirzenius}
\address{Wirzenius: Department of Mathematics and Statistics, Box 68,
FI-00014 University of Helsinki, Finland}
\curraddr{}
\email{henrik.wirzenius@helsinki.fi}
\thanks{}


\begin{abstract}
We exhibit a Banach space $Z$ failing the approximation property, for which there is 
an uncountable family $\mathscr F$ of closed subideals contained in the Banach algebra $\mathcal K(Z)$ 
of the compact operators on $Z$, such that the subideals in $\mathscr F$ are mutually  isomorphic as Banach algebras.  
This contrasts with the behaviour of closed ideals of the algebras $\mathcal L(X)$ of bounded operators on $X$, where 
closed ideals $\mathcal I \neq \mathcal J$ are never isomorphic as Banach algebras. 
We also construct families of non-trivial closed subideals contained in the strictly singular operators $\mathcal S(X)$ for classical spaces such as 
$X = L^p$ with $p \neq 2$, where pairwise isomorphic as well as pairwise non-isomorphic subideals occur.
\end{abstract}

\maketitle

\section{Introduction}\label{intro}
Let $X$ be a Banach space and $\mathcal L(X)$  be the Banach algebra of bounded linear  operators $X \to X$.
It was pointed out in \cite[Added in Proof]{JS21} that 
 \textit{if $\mathcal I$ and $\mathcal J$ are closed ideals of $\mathcal L(X)$ for which there is a
Banach algebra isomorphism $\theta: \mathcal I  \to \mathcal J$, then $\mathcal I = \mathcal J$}.
In other words, distinct closed ideals of 
$\mathcal L(X)$ are never isomorphic as Banach algebras.

We will show that, surprisingly enough, the above property fails in general for closed subideals of $\mathcal L(X)$.
We will adhere to the terminology suggested by Patnaik and Weiss \cite{PW13}, \cite{PW14},  and say  that 
$\mathcal J$ is an $\mathcal I$-subideal of  $\mathcal L(X)$, if $\mathcal J \subset \mathcal I$, where $\mathcal I$ is an ideal of $\mathcal L(X)$
and  $\mathcal J$  is an ideal of $\mathcal I$.  
We are only concerned with \textit{closed linear} subideals, 
that is, $\mathcal J \subset \mathcal I$ are closed linear subspaces of $\mathcal L(X)$,
such that  $US \in \mathcal J$ and $SU \in \mathcal J$ whenever $S \in \mathcal J$ and $U \in \mathcal I$ (and similarly for $\mathcal I \subset \mathcal L(X)$).
It will be convenient to say here that $\mathcal J$ is a \textit{non-trivial} subideal of  $\mathcal L(X)$ if $\mathcal J$ 
is not an ideal of $\mathcal L(X)$. (Note that subideals $\mathcal J$ depend on the intermediary ideal $\mathcal I$, but we will occasionally suppress its role.)
Subideals of $\mathcal L(H)$ for Hilbert spaces $H$ were first considered by Fong and Radjavi in \cite{FR83}.  
In particular, they obtained examples of non-trivial singly generated (but non-closed) $\mathcal K(H)$-subideals $\mathcal J$ of $\mathcal L(H)$,
see e.g. \cite[Theorem 1]{FR83} or \cite[Example 1.3]{PW13}. 

Our main result is based on an example constructed  in \cite[Theorem 4.5]{TW22} for different purposes. 
This produces a family $\{\mathcal I_A:  \emptyset \neq A \varsubsetneq \mathbb N\}$ having the size of the continuum of non-trivial closed $\mathcal K(Z)$-subideals,
for which the subideals $\mathcal I_A$ are mutually isomorphic as Banach algebras. 
Here $\mathcal K(Z)$ denotes the closed ideal of $\mathcal L(Z)$  of the compact operators $Z \to Z$, where the above Banach space
 $Z$ fails to have the approximation property (abbreviated A.P.).
In Section \ref{nonisosub} we obtain, by different methods, families of pairwise non-isomorphic as  well as isomorphic non-trivial closed $\mathcal S(X)$-subideals  
of $\mathcal L(X)$ for classical Banach spaces including $X = L^p(0,1)$, where $1 \le p < \infty$ and $p \neq 2$. 
Above $\mathcal S(X)$ is the closed ideal of the strictly singular operators $X \to X$.
Our results demonstrate that the closed subideals of $\mathcal L(X)$ behave quite differently compared  
with closed ideals.

References \cite{AK06}, \cite{LT77}  and \cite{LT79} will be our standard sources  for undefined concepts related to Banach spaces,
 and \cite{D00} for notions related to Banach algebras.
We use $X \approx Y$ to indicate linearly isomorphic Banach spaces, and $\mathcal A \cong \mathcal B$ for isomorphic Banach algebras (that is, there
is a Banach algebra isomorphism $\theta: \mathcal A \to \mathcal B$). 
Recall that these notions can differ for spaces of operators, as for instance $\mathcal L(L^p)$ and $\mathcal L(\ell^p)$ 
are linearly isomorphic as Banach spaces  for $1 < p < \infty$ and $p \neq 2$ by \cite{AF96}, but they are not 
isomorphic as Banach algebras by Eidelheit's theorem (see below).

\section{Closed subideals of $\mathcal L(X)$ which are isomorphic as Banach algebras}\label{subideals}

Let $X$ and $Y$ be Banach spaces. It is a  classical result of Eidelheit \cite{E40} (see also \cite[Theorem 2.5.7]{D00}) 
that if $\theta: \mathcal L(X) \to \mathcal L(Y)$ is a Banach algebra isomorphism, then there is a linear isomorphism $U \in \mathcal L(X,Y)$ such that
$\theta(S) = USU^{-1}$ for all $S \in \mathcal L(X)$. 
Chernoff \cite[Corollary 3.2]{Ch73} (see also \cite[Section 1.7.15]{P94}) established the following extension: 
\textit{Suppose that $\mathcal A \subset \mathcal L(X)$ and $\mathcal B \subset \mathcal L(Y)$ are subalgebras such that 
the bounded finite rank operators $\mathcal F(X) \subset \mathcal A$ and   $\mathcal F(Y) \subset \mathcal B$. If $\theta: \mathcal A \to \mathcal B$ 
is a bijective algebra homomorphism, then there is a linear isomorphism $U \in \mathcal L(X,Y)$ such that
$\theta(S) = USU^{-1}$ for all $S \in \mathcal A$. }  
As a consequence, if $\mathcal I, \mathcal J \subset \mathcal L(X)$ are closed ideals for which 
there is a Banach algebra isomorphism $\theta:  \mathcal I \to \mathcal J$, then $\mathcal I = \mathcal J$  (cf. also Remarks \ref{isosub}.(ii)).
The purpose of  this section is to exhibit Banach spaces $Z$, where the 
above consequence fails very dramatically within a large class of closed $\mathcal K(Z)$-subideals.

Let $\mathcal A(X,Y) = \overline{\mathcal F(X,Y)}$ denote the class of the approximable operators $X \to Y$, 
where the closure is taken in the uniform operator norm. We note for reference that  
\begin{equation}\label{nontriv}
\mathcal A(X) \subset \mathcal J
\end{equation}
for any non-zero closed $\mathcal I$-subideal  $\mathcal J$ of $\mathcal L(X)$, 
see e.g.  \cite[Theorem 2.5.8.(ii)]{D00} or \cite[Remark 6.1]{PW13}.

We proceed to describe the Banach spaces and the closed subideals  from \cite{TW22}.
Let $(X,Y)$ be a pair of Banach spaces such that 
\begin{equation}\label{pair}
X \textrm{ has the A.P., and }  \mathcal A(X,Y) \varsubsetneq \mathcal K(X,Y).
\end{equation}
We recall that $\mathcal A(X,Y) \varsubsetneq \mathcal K(X,Y)$ for some Banach space $Y$
if and only if the dual space $X^*$ fails the A.P., see e.g.  \cite[Theorem 1.e.5]{LT77}. Moreover, there are spaces  $X$ such that $X$ has the A.P.,
 but $X^*$ fails to have the A.P., see e.g.  \cite[Theorem 1.e.7]{LT77}.  

Fix $1 < p < \infty$. For any pair $(X,Y)$ that satisfies condition \eqref{pair} we consider the direct sum  
\begin{equation}\label{sum}
Z_p := \big(\oplus_{j=0}^\infty X_j\big)_{\ell^p},
\end{equation}
where we put $X_0 = Y$ and $X_j = X$ for $j \ge1$ for unity of notation. 
Bounded operators $S \in \mathcal L(Z_p)$ can be represented as operator matrices $S = (S_{m,n})$ with $S_{m,n} = P_mSJ_n$, where 
$P_m: Z_p \to X_m$ and $J_n: X_n \to Z_p$ are the natural projections and inclusions associated to the component spaces of $Z_p$ for 
$m, n \in \mathbb N \cup \{0\}$. 
For any subset $\emptyset \neq A \varsubsetneq \mathbb N$ define
\begin{equation}\label{id}
\mathcal I_{A} := \{S =  (S_{m,n}) \in \mathcal K(Z_p):  S_{0,0} \in \mathcal A(Y), S_{0,k} \in \mathcal A(X,Y) \textrm{ for all } k \in A\}. 
\end{equation}
It is shown in \cite[Theorem 4.5]{TW22} that the family 
\begin{equation}\label{fam}
\mathscr F := \{\mathcal I_{A}:  \emptyset \neq A \varsubsetneq \mathbb N\}
\end{equation}
has the following properties:

\begin{itemize}
\item[(i)]  $\mathcal I_{A}$ is a closed ideal of $\mathcal K(Z_p)$, and  $\mathcal A(Z_p) \varsubsetneq \mathcal I_{A} \varsubsetneq \mathcal K(Z_p)$ 
for $\emptyset \neq A \varsubsetneq \mathbb N$. 

\item[(ii)]  $\mathcal I_{A}$ is a left ideal of $\mathcal L(Z_p)$ but not a right  ideal  of $\mathcal L(Z_p)$ for $\emptyset \neq A \varsubsetneq \mathbb N$, 
see \cite[Remark 4.8]{TW22} and \cite[Remarks 6.2]{W23}. In particular, $\mathcal I_{A}$ is a non-trivial closed 
$\mathcal K(Z_p)$-subideal of $\mathcal L(Z_p)$ for $\emptyset \neq A \varsubsetneq \mathbb N$.

\item[(iii)]  if $A \subset B$, then $\mathcal I_{B} \subset \mathcal I_{A}$,  and  $\mathcal I_{A} \neq \mathcal I_{B}$  whenever  $A \neq B$. 
\end{itemize}

\smallskip

We stress that above (i)-(iii) hold for the spaces $Z_p$ in \eqref{sum} which are
obtained from any pair $(X,Y)$ that satisfies \eqref{pair}. We will later impose further conditions on $X$ or $Y$, and
in our main result it is assumed that $X$ also satisfies
\begin{equation}\label{extra}
\big(\oplus_{n=1}^\infty X\big)_{\ell^p}  \approx X,
\end{equation}
whence also  $X \oplus X \approx X$. A typical way to achieve this is as follows: if $X_0$ is any space such that $X_0$ has the A.P., but  $X^*_0$ fails the A.P., then 
$\mathcal A(X_0,Y) \varsubsetneq \mathcal K(X_0,Y)$ holds for some space $Y$  by \cite[Theorem 1.e.5]{LT77}.
Let $X =(\oplus_{n=1}^\infty  X_0)_{\ell^p}$. 
It is not difficult to check that $X$ has the A.P. and $\mathcal A(X,Y) \varsubsetneq \mathcal K(X,Y)$, and in addition that
$X  \approx (\oplus_{n=1}^\infty  X)_{\ell^p}$ holds.

Our main result highlights surprising features of the non-trivial closed $\mathcal K(Z_p)$-subideals from the above  family $\mathscr F$. 
This  answers a query of  Gideon Schechtman (private communication). 

\begin{theorem}\label{alliso}
Fix $1 < p < \infty$, and let $Z_p$ be as in \eqref{sum}, where the pair $(X,Y)$ satisfies \eqref{pair} and $X$ satisfies \eqref{extra}.
Then all the non-trivial closed subideals from the family $\mathscr F$ defined by  \eqref{fam}  are mutually isomorphic as Banach algebras, that is,
\[
\mathcal I_A \cong \mathcal I_B \quad \textrm{ for all } \emptyset\neq A, B\subsetneq\mathbb N.
\]
\end{theorem}

Before the proof we comment on the form of Banach algebra isomorphisms between closed subideals. 
Let $X$ be any Banach space and suppose that $\mathcal I, \mathcal J \subset \mathcal L(X)$ are non-zero closed subideals.
Recall from \eqref{nontriv} that non-trivial closed subideals of $\mathcal L(X)$ are closed subalgebras that contain the approximable operators $\mathcal A(X)$.
Hence, if $\theta: \mathcal I \to \mathcal J$ is a Banach algebra isomorphism,
then by  \cite[Corollary 3.2]{Ch73}  there is a linear isomorphism $U \in \mathcal L(X)$, such that the restriction to $\mathcal I$ 
of the inner automorphism 
\[
\psi(S) = USU^{-1}, \quad S \in \mathcal L(X),
\]
equals $\theta: \mathcal I \to \mathcal J$.  In the proof of Theorem  \ref{alliso} we will construct inner automorphisms
$\psi$ of $\mathcal L(Z_p)$ for which $\psi(\mathcal I_A) = \mathcal I_B$. The novel feature is that such a phenomenon is possible among  
non-trivial closed subideals of $\mathcal L(Z_p)$, whereas it is impossible for the smaller class of 
the closed ideals of $\mathcal L(X)$ for any $X$.

\smallskip

The argument will be split into auxiliary steps, where we first construct Banach algebra isomorphisms $\mathcal I_A \cong \mathcal I_B$ for
various basic combinations of the cardinalities  of $A$ and $A^c = \mathbb N \setminus A$, respectively of $B$ and $B^c$. 
In the final step we deduce Theorem \ref{alliso} from these lemmas. Let $\vert A \vert \in\mathbb N\cup\{\infty\}$ denote the cardinality 
of the non-empty set $A \subset \mathbb N$. The spaces $X$ and $Y$ are as in the definition of $Z_p$.

\begin{lemma}\label{perm}
Suppose that $\emptyset\neq A,B\varsubsetneq\mathbb N$  are subsets for which there is a 
permutation $\sigma: \mathbb N \to \mathbb N$ such that $\sigma(A) = B$. Then 
\[
\mathcal I_A\cong \mathcal I_B.
\]
The assumption is satisfied  if (and only if) $|A|=|B|\in\mathbb N\cup\{\infty\}$ and $|A^c|=|B^c|\in\mathbb N\cup\{\infty\}$.
\end{lemma}

\begin{proof}
Define the linear isometry $U \in \mathcal L(Z_p)$ by 
\[
U(y,x_1,x_2,\ldots) = (y,x_{\sigma(1)},x_{\sigma(2)},\ldots), \quad  (y,x_1,x_2,\ldots)  \in Z_p.
\]
Clearly  $U$ is a linear isomorphism $Z_p \to Z_p$, whose inverse $U^{-1}$ satisfies 
\[
U^{-1}(y,x_1,x_2,\ldots) = (y,x_{\sigma^{-1}(1)},x_{\sigma^{-1}(2)},\ldots), \quad  (y,x_1,x_2,\ldots)  \in Z_p.
\]
Let  $\theta(S) = USU^{-1}$ for $S \in \mathcal L(Z_p)$.
It follows that $\theta$ is a Banach algebra isomorphism $\mathcal L(Z_p) \to \mathcal L(Z_p)$, as well as  $\mathcal K(Z_p) \to \mathcal K(Z_p)$.
Its inverse $\theta^{-1}$ on $\mathcal L(Z_p)$ has the form  $\theta^{-1}(T) = U^{-1}TU$ for $T \in \mathcal L(Z_p)$.
It will be enough to verify the following
 
 \smallskip

\underline{Claim}. $\theta(\mathcal I_{A}) \subset \mathcal I_{B}$ and $\theta^{-1}(\mathcal I_{B}) \subset \mathcal I_{A}$. 

\smallskip

Namely, in this event the restriction of $\theta$ to $\mathcal I_{A}$ will be a Banach algebra isomorphism $\mathcal I_{A} \to \mathcal I_{B}$:
for any $T \in \mathcal I_{B}$ one has
$T = \theta (\theta^{-1}(T))$, where $\theta^{-1}(T) \in \mathcal I_{A}$, so that $\theta(\mathcal I_{A}) =  \mathcal I_{B}$.

Towards $\theta(\mathcal I_{A}) \subset \mathcal I_{B}$ we will verify that for any $S \in \mathcal I_{A}$ we have 
$P_0(USU^{-1})J_0 \in \mathcal A(Y)$ and $P_0(USU^{-1})J_r \in \mathcal A(X,Y)$ for any $r \in B$.
Suppose that $y \in Y$ is arbitrary. In this case
\[
SU^{-1}J_0y = SU^{-1}(y,0,0,\ldots) = S(y,0,0,\ldots) = (S_{0,0}y, S_{1,0}y,\ldots),
\]
so that $P_0(USU^{-1})J_0 = S_{0,0} \in \mathcal A(Y)$ since $P_0U =P_0$. 

Next, let $r = \sigma(k) \in B = \sigma(A)$, where $k = \sigma^{-1}(r) \in A$.
If $x_r \in X_r$, then $SU^{-1}J_rx_r = SJ_kx_r$, so that 
\[
SU^{-1}J_rx_r =  (S_{0,k}x_r, S_{1,k}x_r,\ldots).
\]
It follows that  $P_0(USU^{-1})J_r = S_{0,k} \in \mathcal A(X,Y)$,
because $S \in \mathcal I_{A}$ and $k \in A$.

The second inclusion $\theta^{-1}(\mathcal I_{B}) \subset \mathcal I_{A}$ can be deduced from the symmetry. 
Namely,  the inverse permutation $\sigma^{-1}$, for which $\sigma^{-1}(B) = A$, corresponds to the Banach algebra isomorphism  
$\psi(S) =  U^{-1}SU$  for $S \in \mathcal L(Z_p)$. 
The first part of the Claim implies that 
 $\psi(\mathcal I_{B}) \subset \mathcal I_{A}$, where $\psi = \theta^{-1}$. 
  
Finally, if $|A|=|B|\in\mathbb N\cup\{\infty\}$ and $|A^c|=|B^c|\in\mathbb N\cup\{\infty\}$, then there is a permutation 
 $\sigma: \mathbb N \to \mathbb N$ such that  $\sigma(A) = B$ (and $\sigma(A^c) = B^c$).
 \end{proof}
 
 Put $[r]=\{1,\ldots,r\}$ for $r \in \mathbb N$.
 
 \begin{lemma}\label{2}
Suppose that $X\oplus X \approx X$. Then for all $r,s\in\mathbb N$ the following hold:
\begin{enumerate}
\item[(a)] $\mathcal  I_{[r]}\cong \mathcal  I_{[s]}$
\item[(b)] $\mathcal  I_{[r]^c}\cong \mathcal I_{[s]^c}$.
\end{enumerate}
\end{lemma}

\begin{proof}
Let $V:X\to X\oplus X$ be a linear isomorphism.

(a) It will be enough to show that $\mathcal  I_{[r]}\cong  \mathcal  I_{[r+1]}$ for all $r\in\mathbb N$. Namely, if $r <s$, then $\mathcal  I_{[r]}\cong \mathcal  I_{[s]}$
follows by transitivity.

Let $r\in\mathbb N$ and define the  bounded linear isomorphism $U \in \mathcal L(Z_p)$ by 
\[
U(y,x_1,x_2,\ldots)=(y,x_1,\ldots, x_{r-1},V(x_r),x_{r+1},\ldots), \quad  (y,x_1,x_2,\ldots)  \in Z_p,
\]
whose inverse map is 
\[
U^{-1}(y,x_1,x_2,\ldots)=(y,x_1,\ldots, x_{r-1},V^{-1}(x_r,x_{r+1}),x_{r+2},\ldots), \quad  (y,x_1,\ldots)  \in Z_p.
\]
Let $\widetilde{J}_k:X\to X\oplus X$ denote the inclusion maps and $\widetilde{P}_k$ the corresponding projections for $k=1,2$ (relative to $X \oplus X$). 
Observe  that 

\begin{equation}\label{eq1}
UJ_k=\begin{cases}

J_k& \text{if }k\le r-1,\\

J_r\widetilde{P}_1V+J_{r+1}\widetilde{P}_2V&\text{if } k=r,\\

J_{k+1}&\text{if }k>r,
\end{cases}\end{equation}
\begin{equation}\label{eq2}
U^{-1}J_k=\begin{cases}

J_k&\text{if } k\le r-1,\\

J_rV^{-1}\widetilde{J}_1&\text{if } k=r,\\

J_rV^{-1}\widetilde{J}_2&\text{if }k=r+1,\\

J_{k-1}&\text{if }k>r+1.
\end{cases}\end{equation}
Moreover, $P_0U=P_0U^{-1}=P_0$, since the $0$:th component of $Z_p$ is not affected by $U$ or $U^{-1}$.

Let $\theta(S) = USU^{-1}$ for $S \in \mathcal L(Z_p)$, so that $\theta^{-1}(S) = U^{-1}SU$ for $S \in \mathcal L(Z_p)$. It will suffice to verify, as explained
 in the proof of Lemma \ref{perm}, the 
 
 \underline{Claim}. $\theta(\mathcal  I_{[r]})\subset \mathcal  I_{[r+1]}$  and $\theta^{-1}( \mathcal  I_{[r+1]}) \subset \mathcal  I_{[r]}$.

\smallskip

(i) We verify that  $\theta(T) = UTU^{-1}\in \mathcal  I_{[r+1]}$ for any $T\in \mathcal  I_{[r]}$. Note first that $P_0UTU^{-1}J_0=P_0TJ_0\in\mathcal A(Y)$. 
Assume that $k\in[r+1]$. If $k\le r-1$, then 
\[
P_0UTU^{-1}J_k=P_0TJ_k\in\mathcal A(X,Y)
\]
since $U^{-1}J_k = J_k$ by \eqref{eq2}.
If $k=r$, then again by \eqref{eq2} we have
\[P_0UTU^{-1}J_r=P_0TJ_rV^{-1}\widetilde{J}_1\in\mathcal A(X,Y),
\]
since $P_0TJ_r \in \mathcal A(X,Y)$ by assumption.
Finally, if $k=r+1$, then similarly
\[
P_0UTU^{-1}J_{r+1}=P_0TJ_rV^{-1}\widetilde{J}_2\in\mathcal A(X,Y).
\]

(ii) We next verify that  $\theta^{-1}(T) = U^{-1}TU\in \mathcal  I_{[r]}$ for any $T\in \mathcal  I_{[r+1]}$. As above $P_0U^{-1}TUJ_0=P_0TJ_0\in\mathcal A(Y)$. 
Let   $k\in[r]$. If $k\le r-1$, then since $UJ_k = J_k$ by \eqref{eq1} we get that
\[
P_0U^{-1}TUJ_k=P_0TJ_k\in\mathcal A(X,Y)
\]
by assumption. If $k=r$, then from \eqref{eq1} we get that
\[
P_0U^{-1}TUJ_r=P_0T(J_r\widetilde{P}_1V+J_{r+1}\widetilde{P}_2V)\in\mathcal A(X,Y),
\]
since $T\in \mathcal  I_{[r+1]}$ implies that  $P_0TJ_r$ and $P_0TJ_{r+1}$ belong to $\mathcal A(X,Y)$.

\smallskip

(b) Let $U \in \mathcal L(Z_p)$ be the linear isomorphism  from  part (a), and let $\theta(S) = USU^{-1}$ be the corresponding inner automorphism  
$\mathcal L(Z_p) \to  \mathcal L(Z_p)$.  We claim that  also here
\[
\theta(\mathcal I_{[r]^c})\subset \mathcal I_{[r+1]^c} \textrm{ and } \theta^{-1}(\mathcal I_{[r+1]^c}) \subset \mathcal I_{[r]^c}.
\]
As in  part (a) we get that $P_0(\theta(S))J_0 = P_0SJ_0 \in \mathcal A(Y)$ and $P_0(\theta^{-1}(T))J_0=P_0TJ_0\in\mathcal A(Y)$
for any $S \in \mathcal I_{[r]^c}$ and $T \in \mathcal I_{[r+1]^c}$.  

\smallskip

(iii) Suppose  that $S \in \mathcal I_{[r]^c}$ and $s \ge r+2$. From  \eqref{eq2} we have 
\[
P_0(USU^{-1})J_{s} = P_0USJ_{s-1} = P_0SJ_{s-1} \in  \mathcal A(Y)
\]
as $s-1 \ge r+1$ and $S \in \mathcal I_{[r]^c}$.

(iv) Suppose next that $T \in \mathcal I_{[r+1]^c}$ and $s \ge r+1$. From \eqref{eq1} we get that 
\[
P_0(U^{-1}TU)J_s = P_0U^{-1}TJ_{s+1} = P_0TJ_{s+1}  \in  \mathcal A(Y)
\]
as $s+1 \ge r+2$ and $T \in  \mathcal I_{[r+1]^c}$.

This completes the proof of part (b).
\end{proof}

Condition \eqref{extra} on $X$ enables us to find Banach algebra isomorphisms $\mathcal I_A\to \mathcal I_B$ for sets 
$A$ and $B$ of very unequal size. We first isolate two particular cases.

\begin{lemma}\label{3}
Suppose that $X$ satisfies condition \eqref{extra}.  Then the following hold:
\begin{enumerate}
\item[(a)] $\mathcal  I_{\{1\}}\cong \mathcal  I_{\{2,3,4,\ldots\}}$.
\item[(b)] $\mathcal  I_{\{2,3,4,\ldots\}}\cong \mathcal  I_{\{2,4,6,\ldots\}}$.
\end{enumerate}
\end{lemma}

\begin{proof}
Let $V: X \to \big( \oplus_{n=1}^\infty X\big)_{\ell^p}$ be a linear isomorphism.

(a)  We define $U: Z_p \to Z_p$ by 
\[
U(y,x_1,x_2,\ldots) = (y,V^{-1}(x_2,x_3,\ldots),V(x_1)), \quad (y,x_1,x_2,\ldots) \in Z_p,
\]
where $V^{-1}(x_2,x_3,\ldots)$ sits in the first component of $Z_p$. Clearly  $U \in \mathcal L(Z_p)$ is a linear isomorphism for which $U^{-1} = U$.
Let $\psi: \mathcal L(Z_p) \to \mathcal L(Z_p)$ be the Banach algebra isomorphism $\psi(S) = USU$, 
for which $\psi^{-1} = \psi$. Put $B = \{1\}^c$.
 
\underline{Claim}. $\psi(\mathcal I_{\{1\}}) \subset \mathcal I_{B}$ \ and \  $\psi(\mathcal I_{B}) \subset \mathcal I_{\{1\}}$.

(i) Suppose first that $S \in \mathcal  I_{\{1\}}$, so that $S_{0,0}$ and $S_{0,1}$ are approximable operators. Clearly 
$P_0(USU)J_0 = S_{0,0} \in \mathcal A(Y)$. 
Next, let $r \ge 2$ and $x_r \in X_r$ be arbitrary. Then 
\[
SUJ_rx_r = SU(0,0,\ldots,0,x_r,0,\ldots) = S(0,z,0,\ldots) = (S_{0,1}z,S_{1,1}z,\ldots),
\]
where $z = V^{-1}\widetilde{J}_{r-1}x_r$ 
and $\widetilde{J}_k$ is the inclusion $X \to \big( \oplus_{n=1}^\infty  X\big)_{\ell^p}$ into the $k$:th position of the right-hand direct sum.
Deduce that 
\[
 P_0(USU)J_r = S_{0,1}V^{-1}\widetilde{J}_{r-1} \in \mathcal A(X,Y),
 \]
  since $S \in \mathcal  I_{\{1\}}$ and $P_0U = P_0$.
  
(ii)  We next claim that $P_0(UTU)J_1 \in \mathcal A(X,Y)$ for any  $T  \in \mathcal I_{B}$.
For this purpose observe  that 
$\sum_{k=0}^n TJ_kP_k \to T$ as $n \to \infty$ in the operator norm by the proof of \cite[Lemma 4.6]{TW22},
since $T \in \mathcal K(Z_p)$ and $1 < p < \infty$. It follows that 
\[
\Vert \sum_{k=0}^n P_0U(TJ_kP_k)UJ_1 -  P_0UTUJ_1 \Vert \to 0 \textrm{ as } n \to \infty.
\]
By approximation it will suffice to verify that $P_0U(TJ_kP_k)UJ_1 \in \mathcal A(X,Y)$ for all $k \ge 0$, that is, 
$P_0TJ_kP_kUJ_1 \in \mathcal A(X,Y)$ for all $k \ge 0$ (since $P_0U = P_0$). 
Towards this observe that $P_0TJ_k \in \mathcal A(X,Y)$ for $k = 0$ and for $k > 1$ since $T  \in \mathcal I_{B}$.
Moreover, $P_1UJ_1 = 0$ for $k = 1$. Thus $\psi(\mathcal I_{B}) \subset \mathcal I_{\{1\}}$, which completes the proof of part (a).

\smallskip

(b) Define $U: Z_p\to Z_p$ by 
\[
U(y,x_1,x_2,\ldots)=(y,(Vx_1)_1,x_2,(Vx_1)_2,x_3,\ldots), \quad (y,x_1,x_2,\ldots) \in Z_p,
\] 
where $(Vx_1)_k$ denotes the $k$:th component of $Vx_1$ in the direct sum $\big( \oplus_{n=1}^\infty  X\big)_{\ell^p}$.
Then $U \in \mathcal L(Z_p)$ is a linear isomorphism, whose inverse $U^{-1}: Z_p\to Z_p$ is defined by
\[
U^{-1}(y,x_1,x_2,\ldots)=(y,V^{-1}(x_1,x_3,\ldots),x_2,x_4,\ldots), \quad (y,x_1,x_2,\ldots) \in Z_p.
\]

(iii) We first claim that $U^{-1}SU\in \mathcal I_{\{2,3,4,\ldots\}}$ for any  $S\in \mathcal I_{\{2,4,6,\ldots\}}$.
Towards this, note that $P_0U^{-1}=P_0$ and $UJ_0=J_0$. Thus 
\[
P_0U^{-1}SUJ_0=P_0SJ_0\in\mathcal A(Y).
\] 
Suppose next that $r\geq 2$. Observe that  $UJ_r=J_{2r-2}$, and thus 
\[
P_0U^{-1}SUJ_r=P_0SJ_{2r-2}\in\mathcal A(X,Y).\]

(iv) We next claim  that $USU^{-1}\in \mathcal I_{\{2,4,6,\ldots\}}$ for any $S\in \mathcal I_{\{2,3,4,\ldots\}}$.
For this, note again that $P_0U=P_0$ and $U^{-1}J_0=J_0$, so that $P_0USU^{-1}J_0\in\mathcal A(Y)$.
Let $2n\in\{2,4,6,\ldots\}$. In this event $U^{-1}J_{2n}=J_{n+1}$, so that 
\[
P_0 USU^{-1}J_{2n}=P_0SJ_{n+1}\in\mathcal A(X,Y).
\]

Put  $\chi(S) := USU^{-1}$ for $S \in \mathcal L(Z_p)$. By combining parts (iii) and (iv) we deduce 
that $\chi(\mathcal  I_{\{2,3,4,\ldots\}}) = \mathcal  I_{\{2,4,6,\ldots\}}$, so 
$\chi$ yields a Banach algebra isomorphism $\mathcal  I_{\{2,3,4,\ldots\}}\to \mathcal  I_{\{2,4,6,\ldots\}}$.
\end{proof}

We are now in position to  complete the argument of the main result. 

\begin{proof}[Proof of Theorem \ref{alliso}]
By transitivity and symmetry it suffices to show that $\mathcal I_A\cong \mathcal I_{\{1\}}$ for any subset $\emptyset\neq A\varsubsetneq \mathbb N$. 
We consider the cases $|A|<\infty$, $|A^c|<\infty$ and $|A|=|A^c|=\infty$ separately:

\smallskip

\underbar{Case 1}. \textit{Suppose that $|A|=s<\infty$}.
From Lemmas \ref{perm} and \ref{2}.(a) we get that
\[
\mathcal I_A\cong  \mathcal I_{[s]} \cong \mathcal I_{\{1\}}.
\]

\underbar{Case 2}.  \textit{Suppose that  $|A^c|=r<\infty$}.
From Lemmas \ref{perm},  \ref{2}.(b) and \ref{3}.(a) we get that
\[
\mathcal I_A \cong  \mathcal I_{[r]^c} \cong  \mathcal I_{\{1\}^c} \cong  \mathcal I_{\{1\}}.
\]

\underbar{Case 3}.  \textit{Suppose that $|A|=|A^c|=\infty$}.
According to the assumption there is a permutation $\sigma: \mathbb N \to \mathbb N$ such that 
$\sigma(A) = \{2,4,6,\ldots\}$. Hence we find that
\[
\mathcal I_A \cong \mathcal I_{\{2,4,6,\ldots\}}\cong \mathcal I_{\{2,3,4,\ldots\}}\cong \mathcal I_{\{1\}}
\]
from Lemmas \ref{perm}, \ref{3}.(b) and \ref{3}.(a).
\end{proof}

\begin{remarks}\label{isosub}
(i) If $\mathcal K(Z_p)$ is separable in Theorem \ref{alliso}, then it can be verified that $\mathcal K(Z_p)$ has at most continuum many closed subspaces 
(as well as non-trivial $\mathcal K(Z_p)$-subideals). 
Hence the size of the family $\mathscr F$  from \eqref{fam} is as large as possible. 

The pair $(X,Y)$ satisfying \eqref{pair} and \eqref{extra} can be chosen so that $\mathcal K(Z_p)$ is separable. Recall first that 
$\mathcal K(Z_p)$ is separable if and only if $Z_p^*$ is separable, see e.g. \cite[page 272]{TW21}.  
Secondly, if $X^*$ and $Y^*$ are separable, then in \eqref{sum} the dual $Z_p^* = \big(\oplus_{j=0}^\infty X_j^*\big)_{\ell^{p'}}$ is separable.
Here $X_0^* = Y^*$ as well as $X_j^* = X^*$ for $j \ge 1$, and  $p'$ is the dual exponent of $p \in (1,\infty)$.
Next, to choose $X$ we follow the argument of \cite[Theorem 1.e.7.(b)]{LT77}. For this purpose  let $U$ be a separable reflexive space such that $U^*$ fails the A.P.
By \cite[Theorem 1.d.3]{LT77}  there is a Banach space $W$ such that $W^{**}$ has a Schauder basis 
and $W^{**}/W \approx U$. It follows that  $X = W^{**}$ has the A.P., but  $X^* \approx W^* \oplus U^*$ is separable and fails the A.P. 
Apply  \cite[Theorem 1.e.5]{LT77} to pick a Banach space $Y_0$ and $S_0 \in \mathcal K(X,Y_0) \setminus \mathcal A(X,Y_0)$.
By Terzio\u{g}lu's compact factorization theorem \cite{T72} there is a closed subspace $Y \subset c_0$ and a factorization
$S_0 = BS$ with $S \in \mathcal K(X,Y)$. Here $S \notin \mathcal A(X,Y)$ and  $Y^*$ is separable. 

\smallskip

(ii) A  variant of the fact in \cite[Added in Proof]{JS21} implies that the non-trivial subideals $\mathcal I_A \in \mathscr F$ in Theorem \ref{alliso} 
are not isomorphic as Banach algebras to either $\mathcal A(Z_p)$ or $\mathcal K(Z_p)$: 
\textit{Suppose that $X$ is a Banach space, let $\mathcal I$ be a closed ideal of $\mathcal L(X)$ and $\mathcal J$ be a closed subalgebra of  $\mathcal L(X)$
such that $\mathcal A(X) \subset \mathcal J$.
If $\theta: \mathcal I \to \mathcal J$ is a Banach algebra isomorphism, then $\mathcal I = \mathcal J$. 
In particular, if $\mathcal J$ is a non-trivial closed subideal of 
$\mathcal L(X)$, then $\mathcal I$  and $\mathcal J$  are not isomorphic as Banach algebras.}

To see this fact, by  \cite[Corollary 3.2]{Ch73}  there is a linear isomorphism $U \in \mathcal L(X)$ so that $\theta(S) = USU^{-1}$ for $S \in \mathcal I$.
If $T \in \mathcal J$, then there is $S \in \mathcal I$ such that $T = USU^{-1}$, where $USU^{-1} \in \mathcal I$ as $\mathcal I$ is an ideal of 
$\mathcal L(X)$. Thus $\mathcal J \subset \mathcal I$. Conversely, if $S \in \mathcal I$  then $U^{-1}SU \in \mathcal I$, so that 
$S = \theta(U^{-1}SU) \in \mathcal J$. This yields  $\mathcal I = \mathcal J$.
\end{remarks}

Thomas Schlumprecht  asked whether it is possible to identify the closed ideal $[\mathcal I_A]$
of $\mathcal L(Z_p)$ generated by the subideal $\mathcal I_A \in \mathscr F$ for $\emptyset \neq A \varsubsetneq \mathbb N$.
It turns out that $\mathcal I_A$ generate the same closed ideal of $\mathcal L(Z_p)$. We first record another 
general consequence of Chernoff's  result.

\begin{lemma}\label{id1}
Let $X$ be a Banach space and suppose that $\mathcal A \subset \mathcal L(X)$, $\mathcal B \subset \mathcal L(X)$ are closed subalgebras that 
contain $\mathcal F(X)$, for which $\mathcal A \cong \mathcal B$. Then the subalgebras 
$\mathcal A$ and $\mathcal B$ generate the same closed  ideal of $\mathcal L(X)$, that is,
\[
[\mathcal A] = [\mathcal B].
\]
\end{lemma}

\begin{proof}
Let $\theta:  \mathcal A \to \mathcal B$ be a Banach algebra isomorphism. By \cite[Corollary 3.2]{Ch73} there is a linear isomorphism $U \in \mathcal L(X)$ such that 
$\theta(S) = USU^{-1}$ for $S \in  \mathcal A$. If $S \in  \mathcal A$ is arbitrary, then $\theta(S) =  USU^{-1} \in  \mathcal B$, so that $S = U^{-1} \theta(S)U \in [\mathcal B]$.
Deduce that $[\mathcal A] \subset [\mathcal B]$, and by symmetry that $[\mathcal B] \subset [\mathcal A]$.
\end{proof}

 Lemma \ref{id1} together with Theorem \ref{alliso} imply that $[\mathcal I_A] = [\mathcal I_B]$ holds for all non-trivial closed subideals $\mathcal I_A, \mathcal I_B \in \mathscr F$,
 where $\mathscr F$  is given by  \eqref{fam}. 
For this application one requires that  the pair $(X,Y)$ 
of component spaces of  $Z_p$  satisfies \eqref{pair} and that $X$ satisfies \eqref{extra}.
Actually the resulting closed ideal of  $\mathcal L(Z_p)$ can be identified explicitly, and condition \eqref{extra} on $X$ can even be removed.

\begin{prop}\label{fullids}
Suppose that $1 < p < \infty$, and let $Z_p$ be defined by \eqref{sum}, where $(X,Y)$ satisfies \eqref{pair}. Then  
\begin{equation}\label{generate}
[\mathcal I_A]=[\mathcal I]
\end{equation}
for all $\emptyset\neq A\varsubsetneq\mathbb N$, where $\mathcal I:=\{T\in\mathcal K(Z_p)\mid P_0TJ_0\in\mathcal A(Y)\}$. 
\end{prop}

\begin{proof}
Let $\emptyset\neq A\varsubsetneq\mathbb N$ be arbitrary. Since $\mathcal I_A\subset \mathcal I$ it will suffice to verify that $\mathcal I\subset[\mathcal I_A]$. 
Let $T\in \mathcal I$.  Since $T\in\mathcal K(Z_p)$ and $1 < p < \infty$ we know that  
\begin{equation}\label{1}
||\sum_{k=0}^r TJ_kP_k-T||\to 0 \textrm{ as } r\to\infty
\end{equation}
 (see e.g. the proof of \cite[Lemma 4.6]{TW22}). Thus, in order to show that $T\in[\mathcal I_A]$, it will be 
 enough by \eqref{1} to verify that $TJ_kP_k\in[\mathcal I_A]$ for all $k\geq 0$. 
 We need to  consider the following mutually exclusive cases.
 
 \smallskip

\underbar{Case $k=0$}. We know that $P_0(TJ_0P_0)J_0=P_0TJ_0\in\mathcal A(Y)$ by assumption. Moreover, for any $r\in A$ we get that
$P_0(TJ_0P_0)J_r=0$ since $P_0J_r=0$. Thus $TJ_0P_0\in \mathcal I_A$. 

\underbar{Case $k\in A^c$}. Here $TJ_kP_k\in \mathcal I_A$ since $P_0(TJ_kP_k)J_s = 0$ for any $s \in A\cup\{0\}$.

\underbar{Case $k\in A$}. Pick $r\in  A^c$ and let $J_{r,k}:X_r\to X_k$ and $J_{k,r}:X_k\to X_r$ denote the identity operator on $X = X_r = X_k$.
Clearly $J_{r,k}P_rJ_rJ_{k,r}$ is the identity operator $X_k\to X_k$, so that 
\[
TJ_kP_k=(TJ_kJ_{r,k}P_r)(J_rJ_{k,r}P_k).
\]
We claim that $TJ_kJ_{r,k}P_r\in\mathcal I_A$, so that $TJ_kP_k\in[\mathcal I_A]$. In fact, for 
any $s\in A\cup\{0\}$ we have $P_rJ_s=0$, and thus $P_0(TJ_kJ_{r,k}P_r)J_s=0$.
\end{proof}

\begin{remark} 
In Proposition \ref{fullids} there are pairs $(X,Y)$ satisfying \eqref{pair}, for which 
$\mathcal I$ is a non-trivial closed $\mathcal K(Z_p)$-subideal of $\mathcal L(Z_p)$. Hence 
the closed ideal $[\mathcal I]$ of $\mathcal L(Z_p)$ is required on the right-hand side of \eqref{generate}  instead of  $\mathcal I$.

In fact, if $X$ has the A.P. and $X^*$ fails this property, then first apply \cite[Theorem 1.e.5]{LT77}  to pick  $Y_0$ and 
$S_0 \in \mathcal K(X,Y_0) \setminus \mathcal A(X,Y_0)$. 
Let $Y = Y_0 \oplus X$. Then $(X,Y)$ satisfies \eqref{pair}, since $Sx = (S_0x,0)$ for $x \in X$ defines a compact, non-approximable  operator $X \to Y$.
Moreover,  $U(y,x) = x$ for $(y,x) \in Y$ is a non-compact operator $U: Y \to X$ for which 
$SU: Y \to Y$ is compact and non-approximable.
Fix $r \in \mathbb N$ and  define $V \in \mathcal L(Z_p)$ and $T \in \mathcal I$  by $V = J_rUP_0$, respectively $T = J_0SP_r$. Then 
\[
P_0(TV)J_0 = P_0(J_0SP_rJ_rUP_0)J_0 = SU \notin \mathcal A(Y),
\]
that is, $TV \notin \mathcal I$.
\end{remark}

\section{Non-trivial closed $\mathcal S(X)$-subideals}\label{nonisosub}

It is a natural question whether $\mathcal L(X)$ contains non-trivial closed subideals for classical Banach spaces $X$.
Recall that  $\mathcal S(X)$, the class of the strictly singular operators $X \to X$, is a  closed ideal of $\mathcal L(X)$ that satisfies  
$\mathcal K(X) \subset \mathcal S(X)$ for any $X$.  
Here we describe large families of non-trivial closed $\mathcal S(X)$-subideals of  $\mathcal L(X)$ for many classical Banach 
spaces $X$, including $L^p := L^p(0,1)$ with  $p \neq 2$. 
We first briefly discuss closed subideals of Banach algebras.

Let  $\mathcal A$ be a Banach algebra, and suppose that  $\mathcal J \subset \mathcal I \subset \mathcal A$. We say that 
$\mathcal J$ is a closed $\mathcal I$-subideal of $\mathcal A$ if  $\mathcal I$ is a closed ideal of $\mathcal A$ and
$\mathcal J$ is a closed ideal of $\mathcal I$. This setting reveals that the existence 
of non-trivial closed  $\mathcal I$-subideals is related to the absence of approximate identities for $\mathcal I$.
Recall that the net $(e_\alpha) \subset \mathcal I$ is a left approximate identity (LAI)  of
$\mathcal I$ if $y = \lim_{\alpha} e_\alpha y$ for all $y \in \mathcal I$. Right approximate identities (RAI) of $\mathcal I$
are defined analogously.  The following fact  is a variant and reformulation of \cite[Proposition 2.9.4]{D00}.

\begin{lemma}\label{rai}
Suppose that $\mathcal J \subset \mathcal I \subset \mathcal A$, where $\mathcal I$ is a closed ideal of   $\mathcal A$ and 
$\mathcal J$ is a closed $\mathcal I$-subideal of $\mathcal A$.

(i) If  $\mathcal I$ or $\mathcal J$ has a RAI, then $\mathcal J$ is a right ideal of $\mathcal A$.

(ii)  If  $\mathcal I$ or $\mathcal J$  has a LAI, then $\mathcal J$ is a left ideal of $\mathcal A$.
\end{lemma}

\begin{proof}
(i) Suppose that $(e_\alpha)$ is a RAI for $\mathcal I$, and let $x \in \mathcal J$ and $z \in \mathcal A$ be arbitrary. 
It follows that 
\[
xz = \lim_{\alpha} (xe_\alpha)z = \lim_{\alpha} x(e_\alpha z) \in \mathcal J,
\]
since $e_\alpha z \in \mathcal I$ for all $\alpha$ and $\mathcal J$ is a closed ideal of $\mathcal I$. The other cases are similar.
\end{proof}

\begin{remarks}\label{remA}
(i) By Lemma \ref{rai} there are Banach algebras without non-trivial closed subideals. Let 
$\mathcal A$ be a $C^*$-algebra and $\mathcal I \subset \mathcal A$ a closed ideal. It is known that there is a bounded net 
$(e_\alpha) \subset \mathcal I$ which is a LAI as well as a RAI for $\mathcal I$, see e.g. \cite[Proposition 1.8.5]{Di64} or \cite[Theorem 3.2.21]{D00}.
Moreover,  $\mathcal L(X)$ fails to have non-trivial closed subideals by \eqref{nontriv} 
for $X = \ell^p$ with $1 \le p < \infty$ or $X = c_0$, since here $\mathcal K(X)$ is the unique proper closed ideal of $\mathcal L(X)$,
see e.g. \cite[sections 5.1--5.2]{Pi80}.

\smallskip

(ii) If $X$ has the A.P., then $\mathcal K(X) = \mathcal A(X)$, so by \eqref{nontriv} there are no  non-trivial closed 
$\mathcal K(X)$-subideals.
The existence of a LAI or a RAI in $\mathcal K(X)$ is related to the compact approximation property, 
see  \cite[Theorem 2.7]{D86} and \cite[Proposition 7]{Z03}.

\smallskip

(iii) Let $Z_p$ be the Banach space  in \eqref{sum}, where $1 < p < \infty$. 
Lemma \ref{rai} implies that $\mathcal K(Z_p)$ cannot have a RAI. Namely, for any 
$\emptyset \neq A \varsubsetneq \mathbb N$ the closed  $\mathcal K(Z_p)$-subideal  $\mathcal I_{A}$
is not a right ideal of  $\mathcal L(Z_p)$ by property (ii) following \eqref{fam}.
\end{remarks}

For classical spaces $X$ that have the A.P., including $L^p$ for $1 \le p < \infty$ or  $C(0,1)$, 
the algebra  $\mathcal L(X)$ does not contain any non-trivial closed  $\mathcal K(X)$-subideals, see Remarks \ref{remA}.(ii).
The following elementary observation will lead to large families of non-trivial closed $\mathcal S(X)$-subideals of $\mathcal L(X)$, where the 
conditions apply to many classical Banach spaces (see Proposition \ref{elem}). 
 
\begin{prop}\label{2nilsub}
Suppose that $X$ is a Banach space such that
\begin{equation}\label{sing6}
dim\big(\mathcal S(X)/\mathcal K(X)\big) \ge 2,
\end{equation}
\begin{equation}\label{sing4}
UV \in \mathcal K(X) \textrm{ for any } U, V \in \mathcal S(X).
\end{equation}

(i) If  $\mathcal K(X) \varsubsetneq M \varsubsetneq \mathcal S(X)$ is any closed linear subspace, then $M$
is a closed $\mathcal S(X)$-subideal of $\mathcal L(X)$.

(ii) Let $\mathcal K(X) \varsubsetneq M_1, M_2 \varsubsetneq \mathcal S(X)$ be closed linear subspaces. Then the subideals $M_1 \cong M_2$  
if and only if $M_2 = UM_1U^{-1}$ for some linear isomorphism $U \in \mathcal L(X)$.
\end{prop} 

\begin{proof}
(i) If $S \in M$ and $U \in \mathcal S(X)$, then $US \in M$ and $SU \in M$ by  \eqref{sing4}.

\smallskip

(ii) Let  $\theta: M_1 \to M_2$ be a Banach algebra isomorphism. Use \cite[Corollary 3.2]{Ch73}  to find 
a linear isomorphism $U \in \mathcal L(X)$ for which $\theta(S) = USU^{-1} \in M_2$ for $S \in \mathcal M_1$, so that $UM_1U^{-1} \subset  M_2$.
The inner automorphism $\theta: \mathcal L(X) \to \mathcal L(X)$ has the unique inverse $\theta^{-1}$ given by $\theta^{-1}(T) = U^{-1}TU$ for $T \in \mathcal L(X)$.
If $T \in M_2$ is arbitrary, then
\[ 
T= U(U^{-1}TU)U^{-1} = U(\theta^{-1}(T))U^{-1} \in UM_1U^{-1}
\]
as $\theta^{-1}(T) \in M_1$. Conclude that 
$M_2 = UM_1U^{-1}$.

Conversely, suppose that $M_2 = UM_1U^{-1}$ for some linear isomorphism $U \in \mathcal L(X)$. Then the restriction of  the inner automorphism 
$S \mapsto \theta(S) = USU^{-1}$ on $\mathcal L(X)$ is a Banach algebra isomorphism $M_1 \to M_2$.
\end{proof}

Recall that Tarbard \cite{Ta12} constructed a Banach space $X_2$, for  which  
\eqref{sing4} holds and $dim\big(\mathcal S(X_2)/\mathcal K(X_2)\big) = 1$.  
We next provide examples in the above setting of closed linear subspaces $M$, such that  $M$ is not an ideal of $\mathcal L(X)$
and $UMU^{-1} \neq M$ for some linear isomorphism $U$. 
It will be convenient to work on $X \oplus X$, but the spaces $X$ listed below  in Proposition \ref{elem} 
satisfy $X \oplus X \approx X$. 
Note from Eidelheit's theorem that  if  $V: X \oplus X \to X$ is a linear isomorphism, then $S \mapsto \psi(S) = VSV^{-1}$ is a Banach algebra isomorphism 
$\mathcal L(X \oplus X) \to \mathcal L(X)$. Moreover, $\psi(\mathcal K(X \oplus X)) = \mathcal K(X)$ and 
$\psi(\mathcal S(X \oplus X)) = \mathcal S(X)$, so that a non-trivial closed $\mathcal S(X \oplus X)$-subideal $\mathcal I$ transfers to 
a non-trivial closed $\mathcal S(X)$-subideal $\psi(\mathcal I)$ of $\mathcal L(X)$.
We will often write operators $U \in \mathcal L(X\oplus X)$ as 
$U = \begin{bmatrix}
U_{11}&U_{12}\\
U_{21}&U_{22}
\end{bmatrix}$, where $U_{kl} = P_kUJ_l \in \mathcal L(X)$  for $k, l = 1,2$. Here $P_k$ and $J_l$ are the canonical projections and inclusions associated to $X \oplus X$. 
Given closed linear subspaces $M_{ij} \subset \mathcal L(X)$ for $i, j = 1,2$ we denote
\[
\begin{bmatrix}
M_{11}&M_{12}\\
M_{21}&M_{22}
\end{bmatrix} = \Big\{ U = \begin{bmatrix}
U_{11}&U_{12}\\
U_{21}&U_{22}
\end{bmatrix}  \in \mathcal L(X\oplus X):  U_{ij} \in M_{ij} \textrm{ for } i, j = 1,2\Big\}.
\] 
We write  $\mathcal A \ncong \mathcal B$ to indicate non-isomorphic Banach algebras $\mathcal A$ and $\mathcal B$.

\begin{theorem}\label{singularA}
Suppose that $X$ is a Banach space that satisfies \eqref{sing4} and
\begin{equation}\label{sing5}
\mathcal S(X)/\mathcal K(X) \textrm{ is infinite-dimensional}.
\end{equation}

(i) For any closed linear subspace $\mathcal K(X) \varsubsetneq M \varsubsetneq \mathcal S(X)$ put
\begin{equation}\label{elsubid}
\mathcal I(M) = \begin{bmatrix}
M&\mathcal K(X)\\
\mathcal K(X)&\mathcal K(X)
\end{bmatrix}.
\end{equation}
Let $U \in \mathcal L(X \oplus X)$ be the isomorphism
$U(x,y) = (y,x)$ for $(x,y) \in X \oplus X$, and   $\theta(S) = USU^{-1}$ for $S \in \mathcal L(X\oplus X)$.
Then $\mathcal I(M)$ is a non-trivial closed $\mathcal S(X\oplus X)$-subideal of $\mathcal L(X\oplus X)$,
and  $\mathcal J(M) := \theta(\mathcal I(M))$ is also a non-trivial closed $\mathcal S(X\oplus X)$-subideal, for which 
$\mathcal J(M) \cong \mathcal I(M)$ and  $\mathcal J(M) \neq \mathcal I(M)$.

\smallskip

(ii) Let $(T_k) \subset \mathcal S(X)$ be a linearly independent sequence modulo $\mathcal K(X)$, 
and let $M_n \subset \mathcal S(X)$ be the closed linear subspace spanned
by $\{T_j: 1 \le j \le n\} \cup \mathcal K(X)$ for $n \in \mathbb N$. Then $\{\mathcal I(M_n): n \in \mathbb N\}$ is an increasing sequence of 
non-trivial closed $\mathcal S(X\oplus X)$-subideals of $\mathcal L(X\oplus X)$, such that $\mathcal I(M_n) \ncong \mathcal I(M_k)$ 
for any $n \neq k$.

\smallskip

(iii) Suppose that $\mathcal K(X) \varsubsetneq \mathcal J_1, \mathcal J_2 \varsubsetneq \mathcal S(X)$ are closed ideals of $\mathcal L(X)$ such that 
$\mathcal J_1 \neq \mathcal J_2$. Then $\mathcal I(\mathcal J_1)$ and $\mathcal I(\mathcal J_2)$
are non-trivial closed $\mathcal S(X\oplus X)$-subideals 
for which $\mathcal I(\mathcal J_1)  \ncong \mathcal I(\mathcal J_2)$.
\end{theorem}

\begin{proof}
(i) $\mathcal I(M)$ is a closed $\mathcal S(X\oplus X)$-subideal of $\mathcal L(X\oplus X)$
by part (i) of Proposition \ref{2nilsub}.
Fix $T \in M \setminus \mathcal K(X)$, let $S = \begin{bmatrix}
T&0\\
0&0
\end{bmatrix}  \in \mathcal I(M)$ and  $V = \begin{bmatrix}
0&0\\
I_X&0
\end{bmatrix}$. It follows that $VS =  \begin{bmatrix}
0&0\\
T&0
\end{bmatrix} \notin \mathcal I(M)$, so that $\mathcal I(M)$ is a non-trivial subideal.

The above $\theta$ defines an inner automorphism $\theta$ on $\mathcal L(X \oplus X)$, so that 
$\mathcal J(M) = \theta(\mathcal I(M)) \cong \mathcal I(M)$. It is not difficult to check that 
\[
\mathcal J(M) = \begin{bmatrix}
\mathcal K(X)&\mathcal K(X)\\
\mathcal K(X)&M
\end{bmatrix} \neq \mathcal I(M),
\]
since $U = \begin{bmatrix}
0&I_X\\
I_X&0
\end{bmatrix}$. As above $\mathcal J(M)$ is not an ideal of $\mathcal L(X\oplus X)$.

\smallskip

(ii) The linear span $M_n = span(\{T_j: 1 \le j \le n\} \cup \mathcal K(X))$ is a closed linear subspace 
that satisfies  $\mathcal K(X) \varsubsetneq M_n \varsubsetneq \mathcal S(X)$ for all $n \in \mathbb N$.
Hence $\mathcal I(M_n) \varsubsetneq \mathcal I(M_{n+1})$ are non-trivial closed $\mathcal S(X\oplus X)$-subideals of $\mathcal L(X\oplus X)$
for $n \in \mathbb N$ by part (i).

Suppose to the contrary that $k < n$ and $\theta: \mathcal I(M_n) \to \mathcal I(M_k)$ is a Banach algebra isomorphism.
By \cite[Corollary 3.2]{Ch73} pick  a linear isomorphism $U \in \mathcal L(X\oplus X)$,
such that $\theta(T)=UTU^{-1}$ for all $T\in \mathcal I(M_n)$. Since $\theta$ is a Banach algebra isomorphism 
of $\mathcal L(X\oplus X)$ for which $\theta(\mathcal K(X\oplus X)) = \mathcal K(X\oplus X)$,
there is an induced linear isomorphism $\mathcal I(M_n)/\mathcal K(X\oplus X) \to \mathcal I(M_k)/\mathcal K(X\oplus X)$.
This cannot happen since $\mathcal I(M_r)/\mathcal K(X\oplus X)$ is $r$-dimensional for $r \in \mathbb N$.

\smallskip

(iii)  $\mathcal I(\mathcal J_1)$ and $\mathcal I(\mathcal J_2)$
are non-trivial closed subideals of $\mathcal L(X\oplus X)$ by part (i).
Suppose that $\theta: \mathcal I(\mathcal J_2) \to \mathcal I(\mathcal J_1)$
is a Banach algebra isomorphism. 
 We claim that  $\mathcal J_1\subset \mathcal J_2$, so that $\mathcal J_1 = \mathcal J_2$ by symmetry.

By \cite[Corollary 3.2]{Ch73}  there is a linear isomorphism $V \in \mathcal L(X\oplus X)$,
such that $\theta(T)=VTV^{-1}$ for all $T\in \mathcal I(\mathcal J_2)$.
Let $S_0 \in \mathcal J_1$ be arbitrary. Since $S = J_1S_0P_1 \in  \mathcal I(\mathcal J_1)$, there is 
$T \in \mathcal I(\mathcal J_2)$ with $S = VTV^{-1}$. 
Write
$T = \begin{bmatrix}
T_{11}&0\\
0&0
\end{bmatrix} + R$, where $T_{11} \in \mathcal J_2$ and $R \in \mathcal K(X\oplus X)$. 
We get that 
\[
S_0 = P_1SJ_1 = (P_1VJ_1)T_{11}(P_1V^{-1}J_1) + P_1(VRV^{-1})J_1 \in \mathcal J_2,
\]
since $\mathcal J_2$ is an ideal of $\mathcal L(X)$.
\end{proof}

\begin{remarks}\label{isosub1}
(i)  Conditions \eqref{sing4} and \eqref{sing5} imply that $\mathcal S(X)$ has neither a LAI nor  a RAI. 
There are also versions of Proposition \ref{2nilsub} and Theorem  \ref{singularA} for certain other pairs of 
closed ideals $\mathcal J \subset \mathcal I \subset \mathcal L(X)$, but we do not pursue this here.

(ii) No examples of non-trivial closed $\mathcal K(X)$-subideals are available along the line of Theorem \ref{singularA},
since it is unknown whether there is a Banach space $X$ such that 
$\mathcal K(X)/\mathcal A(X)$ is non-zero and $2$-nilpotent, see e.g. \cite[Remarks 5.4.(ii)]{TW22}.
\end{remarks}

We briefly recall some classical Banach spaces $X$ to which Theorem \ref{singularA} applies.
The  spaces  listed here are known to satisfy $X \oplus X \approx X$.

\begin{prop}\label{elem}
Conditions  \eqref{sing4} and  \eqref{sing5}  are satisfied by
 $X = L^p$ with $1 \le p < \infty$ and $p \neq 2$, $C(0,1)$, $\ell^\infty$, as well as  $\ell^p \oplus \ell^q$ and $\ell^p \oplus c_0$ 
with $1 \le p  < q < \infty$.
\end{prop}

\begin{proof}
We first make a preliminary observation towards  \eqref{sing5}. Fix a partition 
$\{A_j: j \in \mathbb N\}$ of $\mathbb N$ into infinite subsets, and  suppose that $1 \le p < q < \infty$.
Put $a = \chi_{A} \in \ell^\infty$ for $A \subset \mathbb N$ and let $D_a: \ell^p \to \ell^q$ be the corresponding diagonal operator. 
Here  $D_a \in \mathcal S(\ell^p, \ell^q)$ by the total incomparability of $\ell^p$ and $\ell^q$, see e.g. \cite[Proposition 2.a.2]{LT77}.

\smallskip

\underbar{Claim.}  The family $\{D_{a_{j}}: j \in \mathbb N\}$ is linearly independent modulo $\mathcal K(\ell^p,\ell^q)$.

\smallskip

\noindent Namely, let  $c_1,\ldots, c_m$ be scalars for some  $m \in \mathbb N$, 
and let $(e_k)$ be the unit vector basis in $\ell^p$. If $c_j  \neq 0$ for some  $j \in \{1,\ldots,m\}$, then $\sum_{r=1}^m c_rD_{a_{r}} \notin \mathcal K(\ell^p,\ell^q)$. 
In fact, for $k \in A_j$ one has
\[
\sum_{r=1}^m c_rD_{a_{r}}e_k = c_je_k,
\]
which  fails to have any norm-convergent subsequences. An analogous claim also holds for $\{D_{a_{j}}: j \in \mathbb N\} \subset \mathcal S(\ell^p,c_0)$,
with a similar argument.

Suppose that  $X = L^p$ with $1 \le p < \infty$ and $p \neq 2$. Recall that if $p \neq 1$, then $\ell^p$ and $\ell^2$ are isomorphic to complemented subspaces
of $L^p$. For $p > 2$ let $P: L^p \to \ell^2$ be a surjection and $J: \ell^p \to L^p$ a linear embedding.
Deduce from the above Claim that   $\{JD_{a_{j}}P: j \in \mathbb N\} \subset \mathcal S(X)$ is linearly independent modulo $\mathcal K(X)$.
For $1 < p < 2$ reverse the roles of $2$ and $p$, and for $p=1$ use the facts that $\ell^1$ is complemented in $L^1$ and $\ell^2$ embeds isomorphically into $L^1$.

Simple modifications yield \eqref{sing5} in the other cases. 
For $X = \ell^p \oplus \ell^q$, where $1 \le p < q < \infty$, it suffices to consider $\{D_{a_{j}}: j \in \mathbb N\} \subset \mathcal S(\ell^p, \ell^q) \subset \mathcal S(\ell^p \oplus \ell^q)$, 
and analogously for $X = \ell^p \oplus c_0$.
Finally, \cite[Proposition 1.3]{R69}  and standard duality imply that $\ell^2$ is a quotient space of both $C(0,1)$ and $\ell^\infty$. 
Since $\ell^p$ embeds isometrically into $C(0,1)$ and $\ell^\infty$ for $p > 2$, we may again proceed as above. 

Milman \cite[Teorema 7]{M70} showed that \eqref{sing4} holds for $X = L^p$ with $1 < p < \infty$ and $p \neq 2$. 
For $p = 1$ it was shown by  Pe{\l}czynski \cite[Theorem II.1]{P65} that 
 $\mathcal S(L^1)$ equals the class of weakly compact operators on $L^1$, so that \eqref{sing4} 
 follows from the Dunford-Pettis property of $L^1$, see e.g. \cite[Theorem 5.4.5.(i)]{AK06}. The same argument also applies to $X = C(0,1)$ and $X = \ell^\infty$
 in view of \cite[Theorem I.1]{P65} and \cite[Theorem 5.4.5.(ii)]{AK06}. Finally,
 the identification of the component ideals of $\mathcal S(\ell^p \oplus \ell^q)$ and $\mathcal S(\ell^p \oplus c_0)$, 
 see e.g. \cite[Theorem 5.3.2]{Pi80},  easily yields
\eqref{sing4}  for  $\ell^p \oplus \ell^q$ and $\ell^p \oplus c_0$. 
\end{proof}

Part (iii) of Theorem \ref{singularA} can be improved e.g. for $L^p$, with  $1 < p < \infty$, $p \neq 2$, by using the 
existence of specific large families of closed ideals of $\mathcal L( L^p)$.

\begin{theorem}\label{singularB}
(i)  Let  $X$ be $L^p$ with $1<p<\infty$ and $p\neq 2$, $\ell^p \oplus \ell^q$ with $1 \le p < q < \infty$ or $\ell^p \oplus c_0$ with $1  < p < \infty$.
 Then there are $2^\mathfrak c$ non-trivial closed $\mathcal S(X)$-subideals of $\mathcal L(X)$ that are pairwise non-isomorphic as Banach algebras.

(ii) For $1<p<\infty$ and $p\neq 2$ there is a family $\{\mathcal I_\alpha : \alpha\in\mathcal C\}$ of the size of the continuum of singly generated non-trivial 
closed $\mathcal S(L^p)$-subideals of $\mathcal L(L^p)$ 
that are pairwise non-isomorphic as Banach algebras, but linearly isomorphic as Banach spaces. 
\end{theorem}

\begin{proof}
(i) Johnson and Schechtman \cite[Remark 4.4]{JS21} proved that there are $2^{\mathfrak c}$ different closed ideals  $\mathcal J$ of $\mathcal L(L^p)$ which satisfy 
$\mathcal K( L^p) \subset \mathcal J \subset  \mathcal S( L^p)$.
The claim then follows from part (iii) of Theorem \ref{singularA}, and the facts that $L^p\approx L^p\oplus L^p$ and that $\mathcal S(L^p)$
has at most   $2^\mathfrak c$ closed subspaces, see e.g. \cite[p. 107]{JS21}. Analogous  
results about closed ideals of $\mathcal L(\ell^p \oplus \ell^q)$ and of  $\mathcal L(\ell^p \oplus c_0)$ were shown by 
Freeman, Schlumprecht and Zs{\'a}k \cite[Corollary 9]{FSZ21} (see also its preceding Remark and the Remark on p. 17 of \cite{FSZ21}).

\smallskip

(ii) Let $\mathcal C$ be a continuum of infinite subsets of $\mathbb N$ such that $|\alpha\cap\beta|<\infty$ for all $\alpha,\beta\in\mathcal C$, and let $1<p<2$. 
By \cite[Remarks 4.3 and 4.4]{JS21}, there is a closed complemented subspace $X\subset L^p$ together with operators 
$U,P,T_\alpha \in \mathcal L(X)$ for $\alpha\in\mathcal C$, such that the following properties hold:
\begin{equation}\label{lp1}
T_\alpha UP\in\mathcal S(X)\setminus\mathcal K(X) \textrm{ for all } \alpha\in\mathcal C,
\end{equation}
\begin{equation}\label{lp2}
T_\beta UP \notin [T_\alpha UP] \textrm{ for any } \alpha,\beta\in\mathcal C, \alpha\neq\beta.
\end{equation}
Here $[T_\alpha UP]$ denotes the closed ideal of $\mathcal L(X)$ generated by $T_\alpha UP$.
Write  $L^p = X \oplus M$ and for $\alpha\in\mathcal C$ consider the operators
\begin{equation}\label{192}
R_\alpha=\begin{bmatrix}
T_\alpha UP&0\\
0&0
\end{bmatrix}\in \mathcal S(L^p)\setminus \mathcal K(L^p), \quad S_\alpha=\begin{bmatrix}
R_\alpha &0\\
0&0
\end{bmatrix}\in\mathcal S(L^p\oplus L^p).
\end{equation}

For $\alpha\in\mathcal C$ consider the closed subspace 
\[
M_{\alpha} = \{\lambda R_\alpha +R\mid\lambda\in\mathbb K \textrm{ and } R\in\mathcal K(L^p)\} \subset \mathcal S(L^p)
\]
and let 
\[
\mathcal I_\alpha := \mathcal I(M_{\alpha}) = \{\lambda S_\alpha+R\mid \lambda\in\mathbb K\text{ and }R\in\mathcal K(L^p\oplus L^p)\}
\] 
be as in \eqref{elsubid}. By part (i) of Theorem \ref{singularA} the family $\{\mathcal I_\alpha :  \alpha\in\mathcal C\}$ consists of non-trivial closed 
$\mathcal S(L^p\oplus L^p)$-subideals of $\mathcal L(L^p\oplus L^p)$.
Here $\mathcal I_\alpha$ is a singly generated ideal of $\mathcal S(L^p\oplus L^p)$ according to \eqref{nontriv}, since $L^p$ has a Schauder basis.
It is clear that $\mathcal I_\alpha\approx \mathcal I_\beta$ for all $\alpha,\beta\in\mathcal C$, as 
 the subideals $\mathcal I_\alpha$ for $\alpha\in\mathcal C$ are one-dimensional extensions of $\mathcal K(L^p\oplus L^p)$. 

It follows that $\mathcal I_\alpha\not\cong\mathcal I_\beta$ by a modification of the argument of Theorem \ref{singularA}.(iii). 
In fact, let $\alpha\neq\beta$ and assume that there is a Banach algebra isomorphism $\theta:\mathcal I_\alpha\to \mathcal I_\beta$. Hence
 there is a linear isomorphism $V\in \mathcal L(L^p\oplus L^p)$ such that $\theta(T)=VTV^{-1}$ for $T\in\mathcal I_\alpha$. 
Pick $\lambda\in\mathbb K$ and $R\in\mathcal K(L^p\oplus L^p)$ such that
\[
S_\beta=\theta(\lambda S_\alpha+R)=V(\lambda S_\alpha +R)V^{-1}=\lambda VJ_0T_\alpha UPP_0V^{-1}+VRV^{-1},
\]
where $J_0$ is the inclusion map $X\to L^p\oplus L^p$ (into the first copy of $L^p$), and $P_0$ is the corresponding projection $L^p\oplus L^p\to X$. 
Deduce that 
\begin{equation}\label{2423}
T_\beta UP=P_0S_\beta J_0= \lambda (P_0VJ_0)T_\alpha UP(P_0V^{-1}J_0)+P_0VRV^{-1}J_0\in [T_\alpha UP], 
\end{equation}
which contradicts \eqref{lp2}.

Let $2<q<\infty$ and $p$ be the dual exponent of $q$. For $\alpha\in\mathcal C$ define
\[
\mathcal J_\alpha=\{\lambda S_\alpha^*+R\mid\lambda\in\mathbb K\text{ and }R\in\mathcal K(L^q\oplus L^q)\},
\]
where $S_\alpha$ is given by \eqref{192}. Observe that  $S_\alpha^*\in\mathcal S(L^q\oplus L^q)\setminus\mathcal K(L^q\oplus L^q)$ by 
\cite[Corollary 2]{Weis} or \cite[p. 19]{M70}, so that $\mathcal J_\alpha\subset \mathcal S(L^q\oplus L^q)$ is a non-trivial closed $\mathcal S(L^q\oplus L^q)$-subideal 
of $\mathcal L(L^q\oplus L^q)$ for all $\alpha\in \mathcal C$. Moreover, $\mathcal J_\alpha\approx\mathcal J_\beta$ for all $\alpha,\beta\in\mathcal C$. 

Finally, suppose that $\alpha\neq\beta$ and assume that there is a Banach algebra isomorphism $\theta:\mathcal J_\alpha\to\mathcal J_\beta$,
where as before $\theta(T)=VTV^{-1}$ for all $T\in\mathcal J_\alpha$ and some linear isomorphism $V\in \mathcal L(L^q\oplus L^q)$.
Pick $\lambda\in\mathbb K$ and $R\in\mathcal K(L^q\oplus L^q)$ such that
$S_\beta^*=V(\lambda S_\alpha^* +R)V^{-1}$.
Deduce from  reflexivity that
$S_\beta=(V^*)^{-1}(\lambda S_\alpha+R^*)V^*$,
which leads to a contradiction as in \eqref{2423}. 
\end{proof}

\begin{remark} 
There are versions of  part (ii) of Theorem \ref{singularB} for $X = L^1$, $C(0,1)$ and $\ell^\infty$. For $L^1$
use the family $\{J_pQ: p \in (2,\infty)\} \subset \mathcal S(L^1)$ from \cite[p. 701]{JPS20}, which has similar properties 
to \eqref{lp1} and \eqref{lp2}. Analogous results for $C(0,1)$ are contained in \cite[Corollary 3.2]{JPS20}, and for 
$L^\infty$  in \cite[Theorem 4.2 and Corollary 4.4]{JPS20}, where $\ell^\infty \approx L^\infty$ by \cite[Theorem 4.3.10]{AK06}. 
The details are left to the interested reader.
\end{remark}

\textit{Acknowledgements.} We are indebted to Gideon Schechtman for a query during the conference Functional Analysis in Lille 2022, 
which motivated Theorem \ref{alliso}.
We are also grateful to Thomas Schlumprecht for a question during IWOTA 2020 (August, 2021) which inspired Proposition \ref{fullids},
and to him and Niels Laustsen for suggestions about the 
largest possible size of the family of closed $\mathcal K(Z_p)$-subideals (see Remarks \ref{isosub}.(i)).
We thank the referee for suggestions that improved the presentation.
Henrik Wirzenius gratefully acknowledges the financial support of the Magnus Ehrnrooth Foundation. 

 \bibliographystyle{amsplain}

\end{document}